\documentclass[11pt]{amsart}
\usepackage{amscd,amssymb,amsxtra}
\usepackage[mathscr]{eucal}  
\usepackage{comment}
\usepackage{enumerate}

\makeatletter
\def\section{\@startsection{section}{1}%
  \z@{4\linespacing\@plus\linespacing}{\linespacing}%
  {\normalfont\scshape\centering}}
 \def\subsection{\@startsection{subsection}{2}%
   \z@{1.25\linespacing\@plus.7\linespacing}{.5\linespacing}%
   {\normalfont\bfseries}}
\makeatother

\theoremstyle{definition}
\newtheorem{example}[equation]{Example}

\theoremstyle{remark}

\newtheorem{remark}[equation]{Remark}

\theoremstyle{plain}
\newtheorem{definition}[equation]{Definition}

\newtheorem{theorem}[equation]{Theorem}

\newtheorem{lemma}[equation]{Lemma}

\numberwithin{equation}{section}

\renewcommand{\:}{\colon}

\newcommand{\CC}{{\mathbb C}}

\DeclareMathOperator{\End}{End}

\DeclareMathOperator{\Hom}{Hom}
\DeclareMathOperator{\id}{id}
\DeclareMathOperator{\Map}{Map}

\DeclareMathOperator{\pt}{pt}

\newcommand{\RR}{{\mathbb R}}

\DeclareMathOperator{\Spin}{Spin}

\newcommand{\ZZ}{{\mathbb Z}}

\newcommand{\chiup}{\raise.5ex\hbox{$\chi$}}
\newcommand{\cir}{S^1}
\DeclareMathOperator{\coker}{coker}

\newcommand{\inv}{^{-1}}
\newcommand{\mstrut}{^{\vphantom{1*\prime y}}}

\DeclareMathOperator{\rank}{rank}
\newcommand{\res}[1]{\negmedspace\bigm|_{#1}}
\newcommand{\temsquare}{\raise3.5pt\hbox{\boxed{ }}}

\newcommand{\zmod}[1]{\ZZ/#1\ZZ}

\newcommand{\zt}{\zmod2}

\usepackage[all,2cell,dvips]{xy}\renewcommand{\cir}{\ensuremath{S^1}}
\UseAllTwocells

\DeclareMathOperator{\Pic}{Pic}
\DeclareMathOperator{\String}{String}
\DeclareMathOperator{\cv}{cv}

\DeclareMathOperator{\sign}{sign}
\newcommand{\gpd}{/\!/}
\newcommand{\Gr}[2]{Gr^+_{#1}(\RR^{#2})}
\newcommand{\Hilb}[1]{\mathcal{H}_{#1}}
\newcommand{\M}[1]{\mathcal{M}_{#1}}
\newcommand{\PGK}{\Pic_{\text{g}}K}

\newcommand{\X}{\mathscr{X}}
\newcommand{\Y}{\mathscr{Y}}
\newcommand{\ab}[1]{A_{#1}}
\newcommand{\bX}{\partial X}
\newcommand{\base}[1]{\mathcal{B}(#1)}
\newcommand{\bg}{\bar{\mathfrak{g}}}
\newcommand{\bordR}{\mathcal{BR}iem}
\newcommand{\bordSO}{\mathcal{BSO}}

\newcommand{\fdot}{\dot f}
\newcommand{\field}[1]{\mathcal{F}_{#1}}
\newcommand{\orgp}{\mathcal{O}(G)}
\newcommand{\ort}[1]{\mathfrak{o}(#1)}
\newcommand{\pgk}{pic_{\text{g}}K}
\newcommand{\sbuniv}{\bar\sigma _{\textnormal{univ}}}
\newcommand{\suniv}{\sigma _{\textnormal{univ}}}
\newcommand{\triv}[1]{\underline{#1}}
\newcommand{\twkr}[1]{\tau ^{(KR)}_{#1}}
\newcommand{\tw}[1]{\tau \mstrut_{#1}}
\newcommand{\unit}{\mathbf{1}}
\newcommand{\ztt}{\ZZ/2}

\begin{document}

\abovedisplayskip18pt plus4.5pt minus9pt
\belowdisplayskip \abovedisplayskip
\abovedisplayshortskip0pt plus4.5pt
\belowdisplayshortskip10.5pt plus4.5pt minus6pt
\baselineskip=15 truept
\marginparwidth=55pt




 \title[Consistent Orientation of Moduli Spaces]{Consistent Orientation of
Moduli Spaces} 
 \author[D. S. Freed]{Daniel S.~Freed}
 \thanks{The work of D.S.F. is supported by NSF grant DMS-0603964}
 \address{Department of Mathematics \\ University of Texas \\ 1 University
Station C1200\\ Austin, TX 78712-0257}
 \email{dafr@math.utexas.edu}
  \author[M. J. Hopkins]{Michael J. Hopkins} 
  \thanks{The work of M.J.H. is supported by NSF grant DMS-0306519}
  \address{Department of Mathematics\\Harvard University\\ 1 Oxford Street \\
 Cambridge, MA 02138}
  \email{mjh@math.harvard.edu} 
 \author[C. Teleman]{Constantin Teleman} 
  \thanks{The work of C.T. is supported by EPSRC GR/S06165/01}
 \address{Department of Mathematics \\ University of California \\ 970 Evans
Hall \#3840 \\ Berkeley, CA 94720-3840}  
 \email{teleman@math.berkeley.edu}
 \dedicatory{For Nigel}
 \date{October 19, 2007}
 \begin{abstract} 
 We give an \emph{a priori} construction of the two-dimensional reduction of
three-dimensional quantum Chern-Simons theory.  This reduction is a
two-dimensional topological quantum field theory and so determines to a
Frobenius ring, which here is the twisted equivariant $K$-theory of a compact
Lie group.  We construct the theory via correspondence diagrams of moduli
spaces, which we ``linearize'' using complex $K$-theory.  A key point in the
construction is to \emph{consistently} orient these moduli spaces to define
pushforwards; the consistent orientation induces twistings of complex
$K$-theory.  The Madsen-Tillmann spectra play a crucial role.
 \end{abstract}
\maketitle

In a series of papers~\cite{FHT1,FHT2,FHT3} we develop the relationship
between positive energy representations of the loop group of a compact Lie
group~$G$ and the twisted equivariant $K$-theory~$K^{\tau +\dim G}_G(G)$.
Here $G$~acts on itself by conjugation.  The loop group representations
depend on a choice of ``level'', and the twisting~$\tau $ is derived from the
level.  For all levels the main theorem is an isomorphism of abelian groups,
and for special transgressed levels it is an isomorphism of \emph{rings}: the
fusion ring of the loop group and $K^{\tau +\dim G}_G(G)$ as a ring.  For
$G$~connected with $\pi _1G$~torsionfree we prove in~\cite[\S4]{FHT1}
and~\cite[\S7]{FHT4} that the ring~$K^{\tau +\dim G}_G(G)$ is a quotient of
the representation ring of~$G$ and we calculate it explicitly.  In these
cases it agrees with the fusion ring of the corresponding centrally extended
loop group.  We also treat $G=SO_3$ in~\cite[(A.10)]{FHT4}.  In this paper we
explicate the multiplication on the twisted equivariant $K$-theory for an
arbitrary compact Lie group~$G$.  We work purely in topology; loop groups do
not appear.  In fact, we construct a \emph{Frobenius ring} structure
on~$K^{\tau +\dim G}_G(G)$.  This is best expressed in the language of
topological quantum field theory: we construct a two-dimensional TQFT over
the integers in which the abelian group attached to the circle is~$K^{\tau
+\dim G}_G(G)$.
 
At first glance the ring structure seems apparent. The multiplication map
$\mu \:G\times G\to G$ induces a pushforward on $K$-theory: the Pontrjagin
product.  But in $K$-cohomology the pushforward is the wrong-way, or
\emph{umkehr}, map.  Thus to define it we must $K$-orient the map~$\mu $.
Furthermore, the twistings must be accounted for in the orientations.
Finally, to ensure associativity we must consistently $K$-orient maps
constructed from~$\mu $ by iterated composition.  For connected and simply
connected groups there is essentially a unique choice, but in general one
must work more.  This orientation problem is neatly formulated in the
language of topological quantum field theory.  Cartesian products of~$G$ then
appear as moduli spaces of flat connections on surfaces, and the maps along
which we push forward are restriction maps of the connections to the
boundary.  What is required, then, is a consistent orientation of these
moduli spaces and restriction maps.  The existence of consistent
orientations, which we prove in Theorem~\ref{thm:7}, is in some sense due to
the Narasimhan-Seshadri theorem which identifies moduli spaces of flat
connections with complex manifolds of stable bundles: complex manifolds carry
a canonical orientation in $K$-theory.  Our proof, though, uses only the much
more simple linear statement that the symbol of the de Rham complex on a
surface is the complexification of the symbol of the Dolbeault complex.  As
we explain in~\S\ref{sec:1}, which serves as a heuristic introduction and
motivation, `consistent orientations on moduli spaces' is the topological
analog of `consistent measures on spaces of fields' in quantum field theory.
The latter is what one would like to construct in the path integral approach
to quantum field theory.
 
Our topological construction, outlined in~\S\ref{sec:3}, proceeds via a
universal orientation (Definition~\ref{thm:2}).  The main observation is that
the problem of consistent orientations is a bordism problem, and the relevant
bordism groups are those constructed by Madsen and Tillmann~\cite{MT} in
their formulation of the Mumford conjectures; see~\cite{MW,GMTW} for proofs
and generalizations.  A universal orientation induces a level
(Definition~\ref{thm:4}).  The map from universal orientations to levels is
an isomorphism for simply connected and connected compact Lie groups~$G$, but
in general it may fail to be injective, surjective, or both.  The theories we
construct are parametrized by universal orientations, not by levels.  It is
interesting to ask whether universal orientations also appear in related
topological and conformal field theories as a refinement of the level.
 
The two-dimensional TQFT we construct here is the dimensional reduction of
three-dimensional Chern-Simons theory, refined to have base ring~$\ZZ$ in
place of~$\CC$.  Our construction is \emph{a priori} in the sense that the
axioms of TQFT---the topological invariance and gluing laws---are deduced
directly from the definition.  By contrast, rigorous constructions of many
other TQFTs, such as the Chern-Simons theory, proceed via \emph{generators
and relations}.  Such constructions are based on general theorems which tell
that these generators and relations generate a TQFT: gluing laws and
topological invariance are satisfied.  One can ask if there is an \emph{a
priori} topological construction of Chern-Simons theory using twisted
$K$-theory.  We do not know of one.  In another direction we can extend TQFTs
to lower dimension, so look for a theory in 0-1-2 dimensions which extends
the 1-2 dimensional theory constructed here.  Again, we do not know if there
is an \emph{a priori} construction of that extended theory.

Section~2 of this paper is an exposition of twistings and orientation,
beginning on familiar ground with densities in differential geometry.
Section~4 briefly considers this TQFT for \emph{families} of 1-~and
2-manifolds.  Our purpose is to highlight an extra twist which occurs: that
theory is ``anomalous''.
 
As far as we know, the problem of consistently orienting moduli spaces first
arises in work of Donaldson~\cite{D,DKr}.  He works with anti-self-dual
connections on a 4-manifold and uses excision in index theory to relate all
of the different moduli spaces.  In both his situation and ours the moduli
spaces in question sit inside infinite dimensional function spaces, and the
virtual tangent bundle to the moduli space extends to a virtual bundle on the
function space.  Thus it suffices to orient over the function space, and this
becomes a universal problem.  Presumably our methods apply to his situation
as well, but we have not worked out the details.

We thank Veronique Godin, Jacob Lurie, Ib Madsen, and Ulrike Tillmann for
enlightening conversations.

It is a pleasure and an honor to dedicate this paper to Nigel Hitchin.  We
greatly admire his mathematical taste, style, and influence.  \emph{<Feliz
cumplea\~nos y que cumplas muchos m\'as!}

  \section{Push-Pull Construction of TQFT}\label{sec:1}

 \subsection*{Quantum Field Theory}\label{subsec:1.1}

The basic structure of an $n$-dimensional Euclidean quantum field theory may
be axiomatized simply.  Let $\bordR_n$ be the bordism category whose objects
are closed oriented $(n-1)$-dimensional Riemannian manifolds.  A morphism
$X\:Y_0\to Y_1$ is a compact oriented $n$-dimensional Riemannian manifold~$X$
together with an orientation-preserving isometry of its boundary to the
disjoint union~$-Y_0\sqcup Y_1$, where $-Y_0$~is the oppositely oriented
manifold.  We term~$Y_0$ the \emph{incoming} boundary and~$Y_1$ the
\emph{outgoing} boundary.  A quantum field theory is a functor
from~$\bordR_n$ to the category of Hilbert spaces and trace class maps.  The
functoriality encodes the gluing law; there is also a symmetric monoidal
structure which encodes the behavior under disjoint unions.  There are many
details and subtleties (see~\cite{S1} in this volume, for example), but our
concern is a simpler topological version.  Thus we replace~$\bordR_n$ by the
bordism category~$\bordSO_n$ of smooth oriented manifolds and consider
orientation-preserving diffeomorphisms in place of isometries.  We define an
$n$-dimensional \emph{topological quantum field theory} (TQFT) to be a
functor from $\bordSO_n$ to the category of complex vector spaces.  The
functor is required to be monoidal: disjoint unions map to tensor products.
The functoriality expresses the usual gluing law and the structure of the
domain category~$\bordSO_n$ encodes the topological invariance.  The example
of interest here has an integral structure: the codomain is the category of
abelian groups rather than complex vector spaces.  The integrality reflects
that the theory is a dimensional reduction; see~\cite{F} for a discussion.
 
Physicists often employ a path integral to construct a quantum field theory.
Here is a cartoon version.  To each manifold~$M$ is attached a space~$\field
M$ of fields and so to a bordism $X\:Y_0\to Y_1$ a correspondence diagram 
  \begin{equation}\label{eq:1}
     \xymatrix@!C{&\field{X}\ar[dl]_s \ar[dr]^t \\ \field{Y_0} && \field{Y_1}} 
  \end{equation}
in which $s,t$~are restriction maps.  The important property of fields is
locality: in the diagram 
  \begin{equation}\label{eq:2}
     \xymatrix@!C{&&\field{X'\circ X}\ar[dl]_r \ar[dr]^{r'} \\
     &\field{X}\ar[dl]_s      \ar[dr]^t&&\field{X'}\ar[dl]_{s'} \ar[dr]^{t'}
     \\ \field{Y_0} && \field{Y_1} && 
     \field{Y_2}}   
  \end{equation}
the space of fields~$\field{X'\circ X}$ on the composition of bordisms
$X\:Y_0\to Y_1$ and $X'\:Y_1\to Y_2$ is the fiber product of the maps~$t,s'$.
Fields are really infinite dimensional \emph{stacks}---for example, in gauge
theories the gauge transformations act as morphisms of fields---and the maps
and fiber products must be understood in that sense.
 
The backdrop for the path integral is measure theory.  If there exist
measures~$\mu\mstrut_X,\mu\mstrut_Y$ on the spaces~$\field X,\field Y$ with
appropriate gluing properties, then one can construct a quantum field theory.
Namely, define the Hilbert space~$\Hilb Y=L^2(\field Y,\mu\mstrut_Y)$ and the
linear map attached to a bordism $X\:Y_0\to Y_1$ as the \emph{push-pull}
  \begin{equation*}
     Z_X=t_*\circ s^*\:\Hilb{Y_0}\longrightarrow \Hilb{Y_1}. 
  \end{equation*}
The pushforward~$t_*$ is integration.  Thus if $f\in
L^2(\field{Y_0},\mu\mstrut_{Y_0})$ and $g\in L^2(\field{Y_1},\mu\mstrut
_{Y_1})$, then
  \begin{equation*}
     \left\langle \bar{g},Z_X(f) \right\rangle_{\Hilb{Y_1}} = \int_{\field
     X}\overline{g\bigl(t(\Phi ) \bigr)}\,f\bigl(s(\Phi ) \bigr)\;d\mu\mstrut 
     _X(\Phi ). 
  \end{equation*}
One usually postulates an \emph{action functional} $S_X\:\field X\to\CC$ and
a measure $\tilde\mu\mstrut_X$ such that
$\mu\mstrut_X=e^{-S_X}\tilde\mu\mstrut_X$ and the action satisfies the gluing
law
  $$ S_{X'\circ X}(\Phi )=S_X\bigl(r(\Phi ) \bigr) + S_{X'}\bigl(r'(\Phi )
     \bigr) $$
in~\eqref{eq:2}.  These measures have not been constructed in most examples
of geometric interest.

 \subsection*{Topological Construction}\label{subsec:1.2}

Our idea is to replace the infinite dimensional stack~$\field X$ by a finite
dimensional stack~$\M X\subset \field X$ of solutions to a first order
partial differential equation and to shift from measure theory to algebraic
topology.  Examples of finite dimensional moduli spaces~$\M X$ in
supersymmetric field theory include anti-self-dual connections in four
dimensions and holomorphic maps in two dimensions.  From the physical point
of view the differential equations are the BPS equations of supersymmetry;
from a mathematical point of view they define the minima of a calculus of
variations functional.  In this paper we consider pure gauge theories.  Fix a
compact Lie group~$G$ and for any manifold~$M$ let $\field M$ denote the
stack of $G$-connections on~$M$.  Define~$\M M$ as the stack of \emph{flat}
$G$-connections on~$M$.  If we choose a set~$\{m_i\}\subset M$ of
``basepoints'', one for each component of~$M$, then $\M M$~is represented by
the product of groupoids $\prod_i\left[ \Hom\bigl(\pi _1(M,m_i),G \bigr)\gpd G
\right] $.  A basic property of flat connections is the gluing law
(see~\eqref{eq:2}).

        \begin{lemma}[]\label{thm:1}
 Suppose $X\:Y_0\to Y_1$ and $X'\:Y_1\to Y_2$ are bordisms of smooth
manifolds.  Then $\M{X'\circ X}$~is the fiber product of 
  \begin{equation*}
     \xymatrix@!C{\M{X}\ar[dr]^t &&\M{X'}\ar[dl]_{s'} \\ &\M{Y}} 
  \end{equation*} 
         \end{lemma}

\noindent Roughly speaking, this says that given flat connections on~$X,X'$
and an isomorphism of their restrictions to~$Y$, one can construct a flat
connection on~$X'\circ X$ and every flat connection on~$X'\circ X$ comes this
way.
 
Replace the infinite dimensional correspondence diagram~\eqref{eq:1} with the
finite dimensional correspondence diagram of flat connections: 
  \begin{equation}\label{eq:6}
     \xymatrix@!C{&\M{X}\ar[dl]_s \ar[dr]^t \\ \M{Y_0} && \M{Y_1}} 
  \end{equation}
Whereas the path integral linearizes~\eqref{eq:1} using measure theory, we
propose instead to linearize~\eqref{eq:6} using algebraic topology.  Let
$E$~be a generalized cohomology theory.  To every closed $(n-1)$-manifold we
assign the abelian group 
  \begin{equation*}
     \ab Y=E^{\bullet }(\M Y). 
  \end{equation*}
To a morphism $X\:Y_0\to Y_1$ we would like to attach a homomorphism
$Z_X\:\ab{Y_0}\to \ab{Y_1}$ defined as the push-pull 
  \begin{equation}\label{eq:8}
     Z_X:=t_*\circ s^*\:E^{\bullet }(\M{Y_0})\longrightarrow E^{\bullet
     }(\M{Y_1}) 
  \end{equation}
in $E$-cohomology.  Whereas the path integral requires \emph{measures}
consistent under gluing to define integration~$t_*$, in our topological
setting we require \emph{orientations} of~$t$ consistent with gluing to
define pushforward~$t_*$.  The \emph{consistency} of orientations under
gluing ensures that \eqref{eq:8}~defines a TQFT which satisfies the gluing
law (functoriality).
 
This, then, is the goal of the paper: we formulate the algebro-topological
home for consistent orientations and study a particular example.  Namely,
specialize to~$n=2$ and require that the 1-manifolds~$Y$ and 2-manifolds~$X$
be oriented.  In other words, the domain category of our TQFT is $\bordSO_2$.
For $Y=\cir$ the moduli stack of flat connections is the global quotient
  \begin{equation*}
     \M Y\cong G\gpd G 
  \end{equation*}
of~$G$ by its adjoint action; the isomorphism is the holonomy of a flat
connection around the circle.  Take the cohomology theory~$E$ to be complex
$K$-theory.  The resulting two-dimensional TQFT on oriented manifolds is the
dimensional reduction of three-dimensional Chern-Simons theory for the
group~$G$.  In this case there is a map from consistent orientations to
``levels'' on~$G$; the level is what is usually used to describe Chern-Simons
theory.  A two-dimensional TQFT on oriented manifolds determines a Frobenius
ring and conversely.  The Frobenius ring constructed here is the Verlinde
ring attached to the loop group of~$G$.  The abelian group~$\ab{\cir}$ is a
\emph{twisted} form of~$K(G\gpd G)=K_G(G)$ and its relation to positive energy
representations of the loop group is developed in~\cite{FHT1,FHT2,FHT3}.  In
this paper we describe a topological construction of the ring structure.

 \subsection*{Remarks}\label{subsec:1.3}

 \begin{itemize} 
 \item Let $X$~be the ``pair of pants'' with the two legs incoming and the
single waist outgoing.  Then restriction to the outgoing boundary is the map
$t\:(G\times G)\gpd G\to G\gpd G$ induced by multiplication $\mu \:G\times
G\to G$.  So $Z_X=t_*\circ s^*$, which defines the ring structure in a
two-dimensional TQFT, is pushforward by multiplication on~$G$.  Therefore, we
do construct the Pontrjagin product on~$K^{\tau + \dim G}_G(G)$---here $\tau
$~is the twisting and there is a degree shift as well---and have implicitly
used an isomorphism of twistings $\mu ^*\tau \to \tau \otimes 1 \,+\,
1\otimes \tau $ which, since the TQFT guarantees an associative product,
satisfies a compatibility condition for triple products.  This isomorphism
and compatibility are embedded in our consistent orientation construction.
 \item We do not use the theorem~\cite{A} which constructs a two-dimensional
TQFT from a Frobenius ring.  Rather, our \emph{a priori} construction
manifestly produces a TQFT which satisfies the gluing law, and we deduce the
Frobenius ring as a derived quantity.
 \item Three-dimensional Chern-Simons theory is defined on a bordism category
of manifolds which carry an extra topological structure.  For oriented
manifolds this extra structure is described as a trivialization of~$p_1$ or
signature, or a certain sort of framing.  (For spin manifolds it is described
as a string structure or, since we are in sufficiently low dimensions, an
ordinary framing.)  The two-dimensional reduction constructed here factors
through the bordism categories of oriented manifolds.
 \item The topological push-pull construction extends to families of bordisms
parametrized by a base manifold~$S$.  A choice of consistent orientation
determines this extension to a theory for families of manifolds, albeit an
``anomalous'' theory; see~\S\ref{sec:4} for a discussion.
 \item The pushforward~$t_*$ is only defined if $t\:\M X\to \M {Y_1}$ is a
\emph{representable} map of stacks, i.e., only if the fibers of~$t$ are
spaces---no automorphisms allowed.  This happens only if each component
of~$X$ has a nonempty outgoing boundary.  Therefore, the push-pull
construction only gives a partial TQFT.  We complete to a full TQFT using the
nondegeneracy of a certain bi-additive form; see~\cite[\S17]{FHT3}.
 \item As mentioned earlier, a standard TQFT is defined over the ring~$\CC$
whereas this theory, being a dimensional reduction of a 3-dimensional theory,
is defined over~$\ZZ$.  It is possible to go further and refine the
push-pull construction to obtain a theory over~$K$, where $K$~is the
$K$-theory ring spectrum.  See~\cite{F} for further discussion.
 \item The theory constructed here has two tiers---it concerns 1- and
2-manifolds---so could be termed a `1-2~theory'.  Extensions to
0-1-2~theories, which have three tiers, are of great interest.  The general
structure of such theories has been much studied recently in various
guises~\cite{MS}, \cite{C}, \cite{HL}.  A theory defined down to points is
completely local, and so ultimately has a simpler structure than less local
theories.  We do not know if the push-pull construction here can be extended
to construct a 0-1-2~theory.
 \end{itemize}

\newpage

  \section{Orientation and Twisting}\label{sec:2}

 \subsection*{Ordinary Cohomology}\label{subsec:2.1}

The first example for a differential geometer is de Rham theory.  Let $M$~be
a smooth manifold and suppose it has a dimension equal to~$n$.  An
\emph{orientation} on~$M$, which is an orientation of the tangent
bundle~$TM$, enables integration
  \begin{equation*}
     \int_{M}\:\Omega _c^n(M)\longrightarrow \RR 
  \end{equation*}
on forms of compact support.  Absent an orientation we can integrate twisted
forms, or densities.  The \emph{twisting} is defined as follows.  For any
real vector space~$V$ of dimension~$r$ let $\base{V}$ denote the
$GL_r\RR$-torsor of bases of~$V$.  There is an associated $\ZZ$-graded real
line~$\ort V$ of functions $f\:\base V\to\RR$ which satisfy $f(b\cdot
A)=\sign\det A\cdot f(b)$ for $b\in \base V,\,A\in GL_r\RR$; the degree
of~$\ort V$ is~$r$.  Applied fiberwise this construction yields a flat
$\ZZ$-graded line bundle $\ort V\to M$ for a real vector bundle $V\to M$.
There is a twisted de Rham complex
  \begin{equation}\label{eq:10}
     0 \longrightarrow  \Omega ^{\ort V - r}(M) \xrightarrow{\;\; d \;\;}
     \Omega ^{\ort V - r+1}(M) \xrightarrow{\;\; d \;\;}  \cdots
     \xrightarrow{\;\; d \;\;}  \Omega ^{\ort V}(M)\longrightarrow  0 
  \end{equation}
where $\Omega ^{\ort V + q}(M)$~is the space of smooth sections of the
\emph{ungraded} vector bundle ${\textstyle\bigwedge} ^{r+q}T^*M\otimes \ort
V$.  The cohomology of~\eqref{eq:10} is the twisted de Rham cohomology
$H^{\ort V + \bullet }_{dR}(M)$.  Let $\ort M=\ort{TM}$.  Then integration is
a map
  \begin{equation}\label{eq:11}
     \int_{M}\: \Omega _c^{\ort M}(M)\longrightarrow \RR. 
  \end{equation}
Notice this formulation-notation works if $M$~has several components of
varying dimension: the degree of~$\ort M$ is then the locally constant
function $\dim M\:M\to Z$.

A similar construction works in integer cohomology.  If $\pi \:V\to M$ is a
real vector bundle over a space~$M$ (which needn't be a manifold) we
define~$\ort V$ as the orientation double cover of~$M$ determined by~$V$ and
introduce a $\ZZ$-grading according to the rank of~$V$.  (Note $\rank
V\:M\to\ZZ$ is a locally constant function.)  There is an $\ort V$-twisted
singular complex analogous to~\eqref{eq:10}: cochains in this complex are
cochains on the double cover which change sign under the deck transformation.
The equivalence class of the twisting~$\ort V$ is
  \begin{equation*}
     \bigl[\ort V \bigr]=\bigl(\rank V,w_1(V)\bigr)\in H^0(M;\ZZ)\times
     H^1(M;\zt),  
  \end{equation*}
where $w_1$~is the Stiefel-Whitney class.  The relationship of the twisting
to integration occurs in the Thom isomorphism.  The Thom class $U\in H^{\pi
^*\ort V}_{\cv}(V)$ lies in twisted cohomology with \emph{compact vertical}
support.  Let $B(V),S(V)$ be the ball and sphere bundles relative to a metric
on~$V$.  The \emph{Thom space}~$M^V$ is the pair~$\bigl(B(V),S(V) \bigr)$ or
equivalently, assuming $M$~is a CW complex, the quotient~$B(V)/S(V)$.  The
composite
  \begin{equation*}
     H^{-\ort V+\bullet }(M)\xrightarrow{\;\;\pi ^*\;\;} H^{-\pi
     ^*\ort V+\bullet }(V) \xrightarrow{\;\;U\;\;}H^{\bullet }(M^V) 
  \end{equation*}
is an isomorphism---the Thom isomorphism---a generalization of the suspension
isomorphism.  If $M$~is a \emph{compact} manifold and $i\:M\hookrightarrow
S^n$ a (Whitney) embedding with normal bundle $\nu \to M$, then the
\emph{Pontrjagin-Thom collapse} is the map $c\:S^n \to M^\nu $ defined by
identifying~$\nu $ with a tubular neighborhood of~$M$ and sending the
complement of~$B(\nu )$ in~$S^n$ to the basepoint of~$M^\nu $.  Integration
is then the composite
  \begin{equation}\label{eq:14}
     H^{-\ort \nu +n}(M)\xrightarrow[\cong ]{\;\text{Thom}\;}H^n(M^\nu
     )\xrightarrow{\;\;c^*\;\;}H^n(S^n)\xrightarrow[\cong
     ]{\text{suspension}}\ZZ. 
  \end{equation}
Twistings obey a Whitney sum formula: there is a natural isomorphism 
  \begin{equation*}
     \ort{V_1\oplus V_2}\xrightarrow{\;\;\cong \;\;}\ort{V_1}+\ort{V_2}. 
  \end{equation*}
Applied to $TM\oplus \nu =\triv n$, where $\triv n$~is the trivial bundle of
rank~$n$, we conclude that integration~\eqref{eq:14} is a map
(compare~\eqref{eq:11}) 
  \begin{equation*}
     H^{\ort M}(M)\longrightarrow \ZZ. 
  \end{equation*}
More generally, if $p\:M\to N$ is a proper map there is a pushforward 
  \begin{equation}\label{eq:17}
     p_*\:H^{\ort p+\bullet }(M)\longrightarrow H^{\bullet }(N), 
  \end{equation}
where $\ort p = \ort M - p^*\ort N$.

 \subsection*{$K$-theory}\label{subsec:2.2}

This discussion applies to any multiplicative cohomology
theory.\footnote{also to a cohomology theory defined by a module over a ring
spectrum.}  The only issue is to determine the twisting of a real vector
bundle in that theory.  For complex $K$-theory there are many possible models
for the twisting~$\tw V$ of a vector bundle $V\to M$.  In the
Donovan-Karoubi~\cite{DK} picture $\tw V$~is represented by the bundle of
complex $\zt$-graded Clifford algebras defined by~$V$.  A bundle of algebras
$A\to M$ of this type is considered trivial if $A=\End(W)$ for a $\zt$-graded
complex vector bundle $W\to M$, i.e., if $A$~is Morita equivalent to the
trivial bundle of algebras~$M\times \CC$.  The equivalence class of~$\tw V$
is
  \begin{equation}\label{eq:18}
     [\tw V]=\bigl(\rank V,w_1(V),W_3(V)\bigr)\in H^0(M;\zt)\times
     H^1(M;\zt)\times H^3(M;\ZZ). 
  \end{equation}
Only torsion classes in~$H^3(M;\ZZ)$ are realized by bundles of finite
dimensional algebras, but we have in mind a larger model which includes
nontorsion classes.  (Such models are developed in~\cite{AS}, \cite{FHT1},
\cite{M} among other works.)  There is a Whitney sum isomorphism
  \begin{equation}\label{eq:19}
     \tw{V_1\oplus V_2}\xrightarrow{\;\;\cong \;\;}\tw{V_1}+\tw{V_2}; 
  \end{equation}
the sum of twistings is realized by the tensor product of algebras.  A
$\text{spin}^c$ structure on~$V$ induces an \emph{orientation}, i.e., a
Morita equivalence
  \begin{equation}\label{eq:23}
     \tw{\triv{\rank V}}\xrightarrow{\;\;\cong \;\;}\tw V .
  \end{equation}
An \emph{$A$-twisted} vector bundle is a vector bundle with an $A$-module
structure; it represents an element of twisted $K$-theory.
 
The Whitney formula~\eqref{eq:19} allows us to attach a twisting to any
virtual real vector bundle: set
  \begin{equation}\label{eq:21}
     \tw{-V}=-\tw{V}.
  \end{equation}
Since the Thom space satisfies the stability condition $X^{V\oplus
\triv{n}}\cong \Sigma ^nX^V$, where `$\Sigma $'~denotes suspension, there is
also a Thom spectrum attached to any virtual vector bundle and a
corresponding Thom isomorphism theorem.  An orientation, which is an
isomorphism as in~\eqref{eq:23}, is equivalently a trivialization of the
twisting attached to the reduced bundle $(V-\triv{\rank V})$.

        \begin{remark}[]\label{thm:10}
 There are also twistings of $K$-theory---indeed of any cohomology
theory---which do not come from vector bundles.
        \end{remark}

Suppose $\tau $~is any twisting on a manifold~$N$.  We can put that extra
twisting into the pushforward in $K$-theory associated to a proper map
$p\:M\to N$ (compare~\eqref{eq:17}):
  \begin{equation}\label{eq:68}
     p_*\:K^{(\tau _p + p^*\tau) +\bullet }(M)\longrightarrow K^{\tau +\bullet
     }(N) .
  \end{equation}
Here $\tau _p$~is the twisting $\tau _p=\tau _M - p^*\tau _N$ of the relative
tangent bundle.  In the next section we encounter a situation in which $\tau
_p+p^*\tau $ is trivialized, and so construct a pushforward from untwisted
$K$-theory to twisted $K$-theory; see~\eqref{eq:44}.

Twistings of $K^{\bullet }(\pt)$ form a symmetric monoidal 2-groupoid; its
classifying space~$\PGK$ is thus an infinite loop space.  The notation: $\Pic
K$~is the classifying space of invertible $K$-modules and $\PGK$~the subspace
classifying certain ``geometric'' invertible $K$-modules including twisted
forms of $K$-theory defined by real vector bundles.  As a space there is a
homotopy equivalence
  \begin{equation}\label{eq:20}
     \PGK \sim K(\zt,0)\times K(\zt,1)\times K(\ZZ,3) 
  \end{equation}
with a product of Eilenberg-MacLane spaces, but the group structure on~$\PGK$
is not a product.  The group of equivalence classes of twistings on~$M$ is
the group of homotopy classes of maps~$[M,\PGK]$, which as a set is the
product of cohomology groups in~\eqref{eq:18}. 
 
Let $\pgk$~denote the spectrum whose 0-space is~$\PGK$ and which is a
Postnikov section of the real $KO$-theory spectrum: the
``$\ztt,\,\ztt,\,0,\,\ZZ$'' bit of the ``$\dots
,\,\ztt,\,\ztt,\,0,\,\ZZ,\,0,\,0,\,0,\,\ZZ,\,\dots $'' song.  Thus the
1-space $B\PGK$ of the spectrum~$\pgk$ is a Postnikov section of~$BO$. Also,
let $ko$~denote the connective $KO$-theory spectrum.  Its 0-space is the
group completion of the classifying space of the symmetric monoidal category
of real vector spaces of finite dimension~\cite{S2}.  The map which attaches
a twisting of~$K^{\bullet }(\pt)$ to a real vector space, say via the
Clifford algebra, induces a spectrum map
  \begin{equation}\label{eq:22}
     \tau \:ko\longrightarrow \pgk. 
  \end{equation}

        \begin{remark}[]\label{thm:11}
 If $M$~is a smooth \emph{stack}, then its tangent space is presented as a
graded vector bundle.  An orientation of~$M$ is then an orientation of this
graded bundle, viewed as a virtual bundle by taking the alternating sum of
the homogeneous terms.  The virtual tangent bundle is classified by a map
$M\to ko$ whose composition with~\eqref{eq:22} gives the induced twisting of
complex $K$-theory.  We apply this in~\S\ref{sec:3} to the moduli space of
flat connections on a fixed oriented 2-manifold.  In that case the virtual
tangent bundle is the index of an elliptic complex and the map $M\to ko$ is
computed from the Atiyah-Singer index theorem~\cite{ASi}.
        \end{remark}

The $\ZZ$ part of~$\PGK$, the Eilenberg-MacLane space~$K(\ZZ,3)$, has an
attractive geometric model: gerbes.  Nigel Hitchin has developed beautiful
applications of gerbes in differential geometry~\cite{H3}.  There is a gerbe
model for~$\PGK$ too---one need only add $\zt$-gradings.

  \section{Universal Orientations and Consistent Orientations}\label{sec:3}

 \subsection*{Overview}\label{subsec:3.5}

In this section, the heart of the paper, we define universal orientations
(Definition~\ref{thm:2}) and prove that they exist (Theorem~\ref{thm:7}).  A
universal orientation simultaneously orients the maps~$t$ in~\eqref{eq:6}
along which we pushforward classes in twisted $K$-theory; see~\eqref{eq:71}
for the precise push-pull maps in the theory.  Universal orientations
form a torsor for an abelian group~\eqref{eq:46}.  A universal orientation
determines a level (Definition~\ref{thm:4}), which is the quantity typically
used to label theories.  The map~\eqref{eq:70} from universal orientations to
levels is not an isomorphism in general.
 
We begin with a closed oriented surface~$X$.  The virtual tangent space to
the stack~$\M X$ of flat $G$-connections on~$X$ is the index of a twisted de
Rham complex~\eqref{eq:24}, and we construct a universal
symbol~\eqref{eq:27}---whence universal index---for this operator.  A
trivialization of the universal twisting~\eqref{eq:72} is a universal
orientation, and it simultaneously orients the moduli stacks~$\M X$ for all
closed oriented~$X$. 
 
For a surface~$X$ with (outgoing) boundary we must orient the restriction map
$t\:\M X\to \M{\bX}$ on flat connections.  It turns out that a universal
orientation does this, simultaneously and coherently for all~$X$, as
expressed in the isomorphism~\eqref{eq:43}.  An important step in the
argument is Lemma~\ref{thm:6}, which uses work of Atiyah-Bott~\cite{AB} to
interpret the universal symbol in terms of standard local boundary conditions
for the de Rham complex.

 \subsection*{Closed surfaces}\label{subsec:3.1}

Fix a compact Lie group~$G$ with Lie algebra~$\mathfrak{g}$.  Let $X$~be a
\emph{closed} oriented 2-manifold and $\M X$~the moduli stack of flat
$G$-connections on~$X$.  A point of~$\M X$ is represented by a flat
connection~$A$ on a principal bundle $P\to X$, and the tangent space to~$\M
X$ at~$A$ by the deformation complex
  \begin{equation}\label{eq:24}
     0 \longrightarrow \Omega ^0_X(\mathfrak{g}_P) \xrightarrow{\;\;d_A\;\;}
     \Omega ^1_X(\mathfrak{g}_P) \xrightarrow{\;\;d_A\;\;} \Omega
     ^2_X(\mathfrak{g}_P) \longrightarrow 0, 
  \end{equation}
the de Rham complex with coefficients in the adjoint bundle associated
to~$P$.  This is an elliptic complex.  Its symbol~$\sigma $ satisfies the
reality condition $\sigma (-\xi )= \overline{\sigma (\xi )}$ for~$\xi \in
TX$, since \eqref{eq:24}~is a complex of real differential
operators~\cite{ASi}.  Recall that the symbol of any complex differential
operator lies in $K^0_{\cv}(TX)\cong K^0(X^{TX})$.  The reality condition
gives a lift $\sigma \in KR^0(X^{iTX})$, where the imaginary tangent
bundle~$iTX$ carries the involution of complex conjugation~\cite{At}.  Bott
periodicity asserts that $V\oplus iV$ is canonically $KR$-oriented for any
real vector bundle~$V$ with no degree shift.  In the language of twistings
of~$KR$ this means
  \begin{equation}\label{eq:26}
     \twkr V + \twkr {iV}=0 .
  \end{equation}
Therefore
  \begin{equation}\label{eq:25}
     KR^0(X^{iTX})\xrightarrow[\cong ]{\;\text{Thom}\;} KR^{-\twkr{iTX}}(X)
     \xrightarrow[\cong ]{\;\;\eqref{eq:26}\;\;} KR^{\twkr{TX}}(X)
     \xrightarrow[\cong ]{\;\;\eqref{eq:21}\;\;} KR^{-\twkr{-TX}}(X)
     \xrightarrow[\cong ]{\;\text{Thom}\;}KO^0(X^{-TX}) 
  \end{equation}
from which we locate the symbol~$\sigma \in KO^0(X^{-TX})$.  Note that by
Atiyah duality $KO^0(X^{-TX})\cong KO_0(X)$ and the $KO$-homology group is
well-known to be the home of real elliptic operators.

Now \eqref{eq:24}~is a \emph{universal} operator: its symbol is constructed
from the exterior algebra of~$TX$ and the adjoint representation of~$G$; it
does not depend on details of the manifold~$X$.  Thus it is pulled back from
a universal symbol.  Let $V_n\to BSO_n$ denote the universal oriented
$n$-plane bundle.  The universal symbol lives on the bundle $V_2\to
BSO_2\times BG$, and by~\eqref{eq:25} we identify it as an element
  \begin{equation}\label{eq:27}
     \suniv\in KO^0(BSO_2^{-V_2}\wedge BG_+). 
  \end{equation}
Here $BG_+$~is the space~$BG$ with a disjoint basepoint adjoined and `$\wedge
$'~is the smash product.  Introduce the notation
  \begin{equation*}
     MTSO_n = BSO_n^{-V_n} 
  \end{equation*}
for this Thom spectrum and so write 
  \begin{equation*}
     \suniv\in KO^0(MTSO_2 \wedge BG_+). 
  \end{equation*}
If $f\:X\to BSO_2\times BG$ is a classifying map for~$TX$ and~$P$, and
$\tilde{f}\:X^{-TX}\to MTSO_2\times BG \to MTSO_2\wedge BG$ the induced map
on Thom spectra, then $\sigma =\tilde{f}^*\suniv$.  It is in this sense that
\eqref{eq:27}~is a \emph{universal} symbol.

        \begin{remark}[]\label{thm:12}
 We digress to explain~$MTSO_n$ in more detail.  Let~$\Gr nN$~be the
Grassmannian of oriented $n$-planes in~$\RR^N$ and
  \begin{equation*}
     0\to V_n\to \triv N\to Q_{N-n}\to 0 
  \end{equation*}
the exact sequence of real vector bundles over~$\Gr nN$ in which $V_n$~is the
universal subbundle and $Q_{N-n}$~the  universal quotient bundle.  Denote the
Thom space of the latter as
  \begin{equation*}
     Z_N := \Gr nN^{Q_{N-n}}. 
  \end{equation*}
Then the suspension~$\Sigma Z_N$ is the Thom space of $Q_{N-n}\oplus \triv
1\to \Gr nN$.  But $Q_{N-n}\oplus \triv 1$ is the pullback of
$Q_{N+1-n}\to\Gr n{N+1}$ under the natural inclusion $\Gr nN\hookrightarrow
\Gr n{N+1}$, and in this manner we produce a map $\Sigma Z_N\to Z_{N+1}$.
Whence the spectrum~$MTSO_n=\{Z_N\}_{N\ge 0}$.  The notation
identifies~$MTSO_n$ as an unstable version of the Thom spectrum~$MSO$ and
also alludes to its appearance in the work of Madsen-Tillmann~\cite{MT}.
There are analogous spectra $MTO_n$, $MT\Spin_n$, $MT\String_n$, etc.  If
$F\:S^N\to Z_N$ is transverse to the $0$-section, then $X:=F\inv
(0\text{-section})$ is an $n$-manifold and the pullback of~$Q_{N-n}-\triv{N}$
is stably isomorphic to~$-TX$.  Thus a map\footnote{i.e., a stable map from
the suspension spectrum of~$S$ to~$MTSO_n$} $S\to MTSO_n$ classifies a map
$M\to S$ of relative dimension~$n$ together with a rank~$n$ bundle $W\to M$
and a stable isomorphism\footnote{Thom bordism theories, such as~$MSO$,
retain the information of the stable normal bundle.  Madsen-Tillmann
theories, such as~$MTSO_n$, track the stable \emph{t}angent bundle, which is
one more justification for the~`T' in the notation.} $T(M/S)\cong W$.  An
important theorem of Madsen-Weiss~\cite{MW} and
Galatius-Madsen-Tillmann-Weiss~\cite{GMTW} identifies~$MTSO_n$ as a bordism
theory of fiber bundles rather than a bordism theory of arbitrary maps.
        \end{remark}
 
If a \emph{smooth manifold}~$M$ parametrizes a family of flat $G$-connections
on~$X$---that is, $P\to M\times X$ is a $G$-bundle with a partial flat
connection along~$X$---then there is a classifying map $M\to\M X$ and the
pullback of the stable tangent bundle of~$\M X$ to~$M$ is the index of the
family of elliptic complexes~\eqref{eq:24}.  Note that if we replace the
adjoint bundle~$\mathfrak{g}_P$ in~\eqref{eq:24} by the trivial bundle of
rank~$\dim G$ then the resulting elliptic complex does not vary over~$M$; its
index is a trivializable bundle.  Hence the reduced stable tangent bundle
to~$\M X$ is computed by the de Rham complex coupled to the reduced adjoint
bundle $\bg\mstrut _{P}=\mathfrak{g}\mstrut _P-\triv{\dim G}$.
 
There is a corresponding reduced universal symbol (compare~\eqref{eq:27})
  \begin{equation}\label{eq:31}
     \sbuniv\in KO^0(MTSO_2\wedge BG). 
  \end{equation}
It induces a universal twisting in $K$-theory and a consistent orientation is
constructed by trivializing this twisting.

        \begin{definition}[]\label{thm:2}
 A \emph{universal orientation} is a null homotopy of the composition   
  \begin{equation}\label{eq:72}
     MTSO_2\wedge BG\xrightarrow{\;\sbuniv\;}ko\xrightarrow{\;\;\tau
     \;\;}\pgk. 
  \end{equation}
 Two universal orientations are said to be equivalent if the null homotopies
are homotopy equivalent.  
        \end{definition}

\noindent
   The set of equivalence classes of universal orientations is a torsor
for the abelian group 
  \begin{equation}\label{eq:46}
     \orgp :=[\Sigma MTSO_2\wedge BG,\pgk].
  \end{equation}
We prove in Theorem~\ref{thm:7} below that universal orientations exist.  In
fact, there is a canonical universal orientation, so the torsor of universal
orientations may be naturally identified with the abelian
group~\eqref{eq:46}.  Definition~\ref{thm:2} is designed to orient the
moduli spaces attached to closed surfaces.  In an equivalent form it leads to
the pushforward maps we need for surfaces with boundary and to twisted
$K$-theory of moduli spaces attached to the boundary; see the discussion
preceding~\eqref{eq:39}. 

Return now to the family of partial $G$-connections on~$P\to M\times X$.  The
bundle $P\to M\times X$ is classified by a map $f\:M\to MTSO_2\times BG$ and
the Atiyah-Singer index theorem~\cite{ASi} implies that the index of the
family of operators~\eqref{eq:24} is~$f^*\suniv$.  Therefore, a universal
orientation pulls back to an orientation of the index of~\eqref{eq:24};
c.f. Remark~\ref{thm:11}.  It follows that a universal orientation
simultaneously orients~$\M X$ for all closed oriented 2-manifolds~$X$.

 \subsection*{Surfaces with boundary}\label{subsec:3.2}

As a preliminary we observe two topological facts about~$MTSO_n$.

        \begin{lemma}[]\label{thm:5}
 \begin{enumerate}[(i)] 
 \item $MTSO_1 \simeq S\inv $, the desuspension of the sphere spectrum.
 \item There is a cofibration 
  \begin{equation}\label{eq:69}
     \xymatrix{\Sigma \inv MTSO_{n-1}\ar@{^{(}->}[r]^-{b} & MTSO_n
      \ar@{>>}[r] & (BSO_n)_+.} 
  \end{equation}
 \end{enumerate}
        \end{lemma}

        \begin{proof}
 For~(i) simply observe
  \begin{equation}\label{eq:33}
     MTSO_1 \simeq BSO_1^{-V_1} \simeq \pt^{-\RR} \simeq S\inv .
  \end{equation}
For~(ii) begin with the cofibration built from the sphere and ball bundles of
the universal bundle:
  \begin{equation}\label{eq:35}
     \xymatrix{S(V_n)_+\ar@{^{(}->}[r] & B(V_n)_+\ar@{>>}[r] &
     \bigl(B(V_n),S(V_n)\bigr)}, 
  \end{equation}
Then identify~$BSO_{n-1}$ as the unit sphere bundle~$S(V_n)$ and
write~\eqref{eq:35} in terms of Thom spaces:
  \begin{equation}\label{eq:36}
     \xymatrix{ BSO_{n-1}^{\triv0}\ar@{^{(}->}[r] & BSO_n^{\triv0}\ar@{>>}[r] &
     BSO_n^{V_n}}. 
  \end{equation}
Here $\triv{0}$~is the vector bundle of rank zero.  Now add~$-V_n$ to each of
the vector bundles in~\eqref{eq:36} and note that the restriction of~$V_n$
to~$BSO_{n-1}$ is $V_{n-1}\oplus \triv1$.
        \end{proof}

Consider the diagram 
  \begin{equation}\label{eq:38}
     \xymatrix{\Sigma \inv MTSO_1\wedge BG \ar[r]^-{b} & MTSO_2\wedge BG
     \ar[r]^-{q}\ar[dr]^-{\sbuniv} & (MTSO_2,\Sigma \inv MTSO_1)\wedge
     BG\ar@{.>}[d]^-{\sbuniv'} \\ 
     &&ko}
  \end{equation}
The top row is a cofibration.  From~\eqref{eq:69} we can
replace~$(MTSO_2,\Sigma \inv MTSO_1)$ with~$(BSO_2)_+$.

        \begin{lemma}[]\label{thm:3}
 Define $\sbuniv'$ in~\eqref{eq:38} to be the map $(BSO_2)_+\wedge BG\to ko$
induced by the trivial representation of~$SO_2$ smash with the reduced
adjoint representation of~$G$.  Then the triangle in~\eqref{eq:38} commutes
and the diagram gives a canonical null homotopy of the
composite~$\sbuniv\circ b$.
        \end{lemma}

        \begin{proof}
  Recalling the isomorphisms in~\eqref{eq:25}, and using the fact that the
universal symbol~$\sbuniv$ is canonically associated to a representation
of~$SO_2\times G$, we locate $\sbuniv\in KR^0_{SO_2\times G}(-\RR^2)\mstrut
_c\cong KR^0_{SO_2\times G}(i\RR^2)\mstrut _c$.  (Recall that the involution
on $\CC^2\cong \RR^2\oplus i\RR^2$ is complex conjugation and the
subscript~`c' denotes compact support.)  Similarly, $\sbuniv'\in
KR^0_{SO_2\times G}(\pt)\cong KR^0_{SO_2}(\CC^2)\mstrut _c$.  Using the Thom
isomorphism we identify~$\sbuniv'$ as the difference of classes represented
by ${\textstyle\bigwedge} ^{\bullet }\CC^2\otimes \mathfrak{g}\mstrut _{\CC}$
and ${\textstyle\bigwedge} ^{\bullet }\CC^2\otimes \CC^{\dim G}$, where in
both summands $\theta \in \CC^2$~acts as $\epsilon (\theta )\otimes \id$.
Exterior multiplication~$\epsilon (\theta )$ is exact for~$\theta \not= 0$,
so this does represent a class with compact support.  Also, as $\epsilon $
commutes with complex conjugation it is Real in the sense of~\cite{At}.  It
remains to observe that its restriction under $KR^0_{SO_2\times
G}(\CC^2)\mstrut _c\to KR^0_{SO_2\times G}(i\RR^2)\mstrut _c$ is the
universal symbol~$\sbuniv$ of the de Rham complex coupled to the reduced
adjoint bundle.
        \end{proof}

For any~$n$ a point of the 0-space of the pair of spectra $(MTSO_n,\Sigma
\inv MTSO_{n-1})$ is represented by a map from~$(B^N,S^{N-1})$ into the pair
of Thom spaces attached to
  \begin{equation*}
     \xymatrix{Q_{N-n}\ar[r]\ar[d] & Q_{N-n}\ar[d]\\
     \Gr{n-1}{N-1}\ar@{^{(}->}[r] & \Gr nN} 
  \end{equation*}
for $N$~sufficiently large.  A map of~$(B^N,S^{N-1})$ into this pair which is
transverse to the $0$-section gives a compact oriented $n$-manifold~$M$ with
boundary embedded in~$(B^N,S^{N-1})$, a rank~$n$ bundle $W\to M$ equipped
with a stable isomorphism $TM\cong W$, and a splitting of~$W\res{\partial
M}$ as the direct sum of a rank~$(n-1)$ bundle and a trivial line bundle.
The composition with the boundary map
  \begin{equation}\label{eq:52}
     r\:(MTSO_n,\Sigma \inv MTSO_{n-1})\to MTSO_{n-1} 
  \end{equation}
is the restriction of this data to~$\partial M$.
 
Now let $M$~be a smooth manifold, $X$~a compact oriented 2-manifold with
boundary, and $P\to M\times X$ a principal $G$-bundle with partial flat
connection along~$X$.  This data induces a classifying map $M\to\M X$ and,
forgetting the connection, a classifying map
  \begin{equation}\label{eq:41}
     f\:M\longrightarrow (MTSO_2,\Sigma \inv MTSO_1)\wedge  BG. 
  \end{equation}
There are induced classifying maps $M\to\M{\bX}$ and 
  \begin{equation*}
     \fdot\:M\longrightarrow MTSO_1\wedge  BG 
  \end{equation*}
for the boundary data; here $\fdot=r\circ f$.  View~$X$ as a bordism
$X\:\emptyset \to \bX$; later we incorporate incoming boundary components.
The following key result relates the universal topology above to surfaces
with boundary.

        \begin{lemma}[]\label{thm:6}
  The map $\sbuniv'\circ f\:M\to ko$~classifies the reduced tangent bundle of
the restriction map $t\:\M X\to\M{\bX} $---the bundle $T\M X-t^*T\M{\bX}$ reduced to rank
zero---pulled back to~$M$.
        \end{lemma}

        \begin{proof}
 At a point $A\in M$ there is a flat connection on $P\res{\{A\}\times X}$.
The tangent space to~$\M X$ at that point is computed by the twisted de Rham
complex~\eqref{eq:24}, so is represented by the twisted de Rham
cohomology~$H^{\bullet }_A(X)$.  Similarly, the tangent space to~$\M{\bX}$ at
the restriction~$t(A)$ of~$A$ to the boundary is~$H^{\bullet }_{t(A)}(\bX)$.
From the long exact sequence of the pair~$(X,\bX)$ we deduce that the
difference $T\M X - t^*T\M{\bX}$ at~$A$ is the twisted relative de Rham
cohomology~$H^{\bullet }_A(X,\bX)$.

Now the twisted relative de Rham cohomology is the index of the deformation
complex~\eqref{eq:24} with \emph{relative boundary
conditions}~\cite[\S4.1]{G}.  In other words, we consider the subcomplex of
forms which vanish when restricted to ~$\bX$.  This is an example of a
\emph{local} elliptic boundary value problem.  Atiyah and Bott~\cite{AB}
interpret local boundary conditions in $K$-theory and prove an index
formula.  More precisely, the triple $\bigl(B(TX),\partial B(TX),S(TX)
\bigr)$ leads to the exact sequence 
  \begin{multline*}
     KR\inv \bigl(B(TX)\res{\bX},S(TX)\res{\bX} \bigr)\longrightarrow
     KR^0\bigl(B(TX),\partial B(TX) \bigr) \\[6pt]
      \longrightarrow
     KR^0\bigl(B(TX),S(TX) \bigr) \longrightarrow KR^0
     \bigl(B(TX)\res{\bX},S(TX)\res{\bX} \bigr) 
  \end{multline*}
The symbol~$\sigma $ of an elliptic operator lives in the third group, and
Atiyah-Bott construct a lift to the second group from a local boundary
condition.  (The image of~$\sigma $ in the last group is an obstruction to
the existence of local boundary conditions; the image of the first group in
the second measures differences of boundary conditions.)  The relative
boundary conditions on the twisted de Rham complex are \emph{universal}, so
the corresponding lift of the symbol occurs in the universal setting.  Recall
from the proof of Lemma~\ref{thm:3} the exact sequence~\eqref{eq:38}, now
extended one step to the left: 
  \begin{equation*}
     KR^0_G(i\RR)\mstrut _c \longrightarrow KR^0_{SO_2\times G}(\CC^2)\mstrut
     _c\longrightarrow KR^0_{SO_2\times G}(i\RR^2)\mstrut _c \longrightarrow
     KR_G^0(i\RR^2)\mstrut _c 
  \end{equation*}
The group~$G$ acts trivially in all cases.  The Atiyah-Bott procedure applied
to the relative boundary conditions gives a lift of~$\sbuniv\in
KR^0_{SO_2\times G}(i\RR^2)\mstrut )_c$ to~$KR^0_{SO_2\times G}(\CC^2)\mstrut
_c$.  Recall that $\sbuniv'$, constructed in the proof of Lemma~\ref{thm:3},
is also a lift of~$\sbuniv$.  But by periodicity we find $KR_G^0(i\RR)\mstrut
_c\cong KR^0_G(-\RR)\mstrut _c\cong KO^0_G(-\RR)_c\cong KO^1_G(\pt)$ which
vanishes by~\cite{An}.  Thus the lift of~$\sbuniv$ is unique and $\sbuniv'$
computes the relative twisted de Rham cohomology.  This completes the proof.
        \end{proof}

 \subsection*{The Level}\label{subsec:3.3}

A universal orientation induces a level, which is commonly used to identify
the theory.  One observation arising from this study is that the level is a
derived quantity, and it is the universal orientation which determines the
theory.  We explain the relationship, and deduce the existence of universal
orientations, in this subsection.

To begin we recast Definition~\ref{thm:2} of a universal orientation in a
form suited for surfaces with boundary.  Consider the diagram
  \begin{equation}\label{eq:39}
     \xymatrix{MTSO_2\wedge BG \ar[r]^-q\ar[dr]^-{\sbuniv} & (MTSO_2,\Sigma
     \inv MTSO_1)\wedge BG \ar[d]^{\sbuniv'}\ar[r]^-r & MTSO_1\wedge BG
     \ar@{.>}[d]^-{-\lambda }\\ & ko \ar[r]^-{\tau } & \pgk} 
  \end{equation}
The top row is a cofibration, the continuation of the top row
of~\eqref{eq:38} in the Puppe sequence.  Recall that a universal orientation
is a null homotopy of $\tau \circ \sbuniv=\tau \circ \sbuniv'\circ q$.

        \begin{lemma}[]\label{thm:9}
 A universal orientation is equivalent to a map~$-\lambda $ in~\eqref{eq:39}
and a homotopy from $\tau \circ \sbuniv'$ to~$-\lambda \circ r$.  
        \end{lemma}

\noindent
 The proof is immediate.  In view of~\eqref{eq:33} and adjunction we can
write $\lambda \:\Sigma ^{\infty}BG\to \Sigma \pgk$ as a map from the
suspension spectrum of~$BG$ to the spectrum~$\pgk$, or equivalently as a map
  \begin{equation}\label{eq:40}
     \lambda \:BG\longrightarrow B\PGK 
  \end{equation}
on spaces. 

        \begin{definition}[]\label{thm:4}
 The map~\eqref{eq:40} is the \emph{level} induced by a universal
orientation. 
        \end{definition}

\noindent
 There is a map $K(\ZZ,4)\to B\PGK$ (see~\eqref{eq:20}) and in some important
cases, for example if $G$~is connected and simply connected, the level
factors through~$K(\ZZ,4)$, i.e., the level is a class~$\lambda \in H^4(BG)$.
 
A universal orientation is more than a level: it is a map $-\lambda
\:MTSO_1\wedge BG\to\pgk$ together with a homotopy of~$-\lambda \circ r$
and~$\tau \circ \sbuniv'$ in~\eqref{eq:39}.  Our next result proves that
universal orientations exist.

        \begin{theorem}[]\label{thm:7}
 There is a canonical universal orientation~$\mu $.  The corresponding
level~$h$ is the negative of $BG\to BO\to B\PGK$, where the first map is
induced from the reduced adjoint representation~$\bg$ and the second is
projection to a Postnikov section.
        \end{theorem}

        \begin{proof}
 Since complex vector spaces have a canonical $K$-theory orientation, the
composite map $\xymatrix{k\ar[r] & ko \ar[r]^-{\tau } &\pgk}$ is null, where
$k$~is the connective $K$-theory spectrum.  (See the text
preceding~\eqref{eq:22}.)  Therefore, a universal orientation is given by
filling in the left dotted arrow in the diagram
  \begin{equation}\label{eq:47}
     \xymatrix{MTSO_2\wedge BG \ar[r]^-q\ar@{.>}[d] &
     (MTSO_2,\Sigma \inv MTSO_1)\wedge BG \ar[d]^{\sbuniv'}\ar[r]^-r &
     MTSO_1\wedge BG \ar@{.>}[d]^-{-\lambda }\\ k\ar[r]& ko \ar[r]^-{\tau } &
     \pgk} 
  \end{equation}
and specifying a homotopy which makes the left square commute.  There is a
natural choice: smash the $K$-theory Thom class $U\:MTU_1\simeq MTSO_2\to k$
with the complexified reduced adjoint representation~$\bg\mstrut _{\CC}$.
This is the universal rewriting of de Rham on a Riemann surface in terms of
Dolbeault, at least on the symbolic level.  In terms of the proof of
Lemma~\ref{thm:3}, the map~$\sbuniv$, restricted to~$MTSO_2$, is the exterior
algebra complex~$({\textstyle\bigwedge} ^{\bullet }\CC^2,\epsilon )$
in~$KR^0_{SO_2}(i\RR^2)_c$.  Write $\RR^2=L_{\RR}$ for the complex
line~$L=\CC$ and $\CC^2\cong \RR^2\otimes \CC\cong L\oplus
\overline{L}$. Then the symbol complex at~$\theta \in i\RR^2$,
  \begin{equation*}
     \xymatrix@!{\CC \ar[r]^-{\epsilon (\theta )} &L\oplus \overline{L}
     \ar[r]^-{\epsilon (\theta )} & L\otimes \overline{L}}, 
  \end{equation*}
is the realification of the complex 
  \begin{equation*}
     \xymatrix@!{\CC \ar[r]^-{\epsilon (\theta )} & L}
  \end{equation*}
which defines the $K$-theory Thom class.  Tensor with the complexified
reduced adjoint representation~$\bg\mstrut _{\CC}$ to complete the argument.

To compute the level of~$\mu $ we factorize~$\tau $ as $\xymatrix{ko
\ar[r]^-{\eta }&\Sigma \inv ko \ar[r] &\pgk}$, where the first map is
multiplication by $\eta \:S^0\to S\inv $ and the second is projection to a
Postnikov section.  The map~$\eta $ fits into the Bott sequence $k\to ko\to
\Sigma \inv ko$, and so we extend~\eqref{eq:47} to the left:
  \begin{equation}\label{eq:51}
     \xymatrix{&\Sigma \inv MTSO_1\wedge BG\ar[r] \ar@{.>}[d]^{\alpha }
     &MTSO_2\wedge BG 
     \ar[r]^-q\ar[d]^{U\wedge \bg\mstrut _{\CC}} & (MTSO_2,\Sigma \inv 
     MTSO_1)\wedge BG \ar[d]^{\sbuniv'}\\  
     \Sigma \inv ko\ar[r] &\Sigma ^{-2}ko\ar[r] &k\ar[r]& ko } 
  \end{equation}
The homotopy which expresses the commutativity of the right hand square
induces the map~$\alpha $ in this diagram, and the map~$-\lambda $ induced
in~\eqref{eq:47} is the suspension of~$\alpha $.  We claim that there is a
unique~$\alpha $, up to homotopy, which makes the left square
in~\eqref{eq:51} commute.  For the difference of any two choices for~$\alpha
$ is a map $\Sigma \inv MTSO_1\wedge BG\to\Sigma \inv ko$, and the homotopy
classes of such maps form the group~$KO^1(BG)$ which vanishes~\cite{An}.  It
is easy to find a map~$\alpha $ as follows.  Since (Lemma~\ref{thm:5}(i))
$\Sigma \inv MTSO_1\simeq S^{-2}$, the upper left map is the inclusion of the
bottom cell of $MTSO_2\wedge BG$.  The composite $\Sigma \inv MTSO_1\wedge
BG\simeq \Sigma ^{-2}BG\to k$ factors as $\Sigma ^{-2}BG\to\Sigma ^{-2}k\to
k$, where the first map is the double desuspension of~$\bg\mstrut _{\CC}$ and
the second Bott periodicity.  Choose~$\alpha $ to be the double desuspension
of~$\bg$, the \emph{real} reduced adjoint representation.  This completes the
proof. 
        \end{proof}

Since equivalence classes of universal orientations form a torsor for the
group~$\orgp$ in~\eqref{eq:46}, the canonical universal orientation
identifies the torsor of universal orientations with~$\orgp$.  Notice the
natural map  
  \begin{equation}\label{eq:70} 
     \ell \:\orgp=[\Sigma MTSO_2\wedge BG,\pgk]\longrightarrow [MTSO_1\wedge
     BG,\pgk]\cong [BG,B\PGK] 
  \end{equation}
from universal orientations to levels.  If~$g\in \orgp$, then the level
of~$\mu + g $ is~$\ell (g)-h$.  If $G$~is connected, simply connected,
and simple, then $[BG,B\PGK]\cong H^4(BG;\ZZ)\cong \ZZ$ and $h$~is the dual
Coxeter number of~$G$ times a generator.  Then $g\mapsto \ell (g)- h$ is a
version of the ubiquitous ``adjoint shift''. 

        \begin{remark}[]\label{thm:13}
 For any~$G$ the top homotopy group of~$\Map(\Sigma MTSO_2,\pgk)$ and of
~$B\PGK$ is~$\pi _4$, which is infinite cyclic.  So there is a homomorphism
of~$H^4(BG;\ZZ)$ into the domain and codomain of~\eqref{eq:70}, and on these
subspaces $\ell$~is an isomorphism.  This means that we can \emph{change} a
consistent orientation by an element of~$H^4(BG;\ZZ)$, and the level changes
by the same amount.
        \end{remark}

 \subsection*{The pushforward maps}\label{subsec:3.4}

Suppose we have chosen a universal orientation with level~$\lambda $.  Let
$X$~be a compact oriented 2-manifold with boundary.  We can work as before
with a family of flat connections on~$X$ parametrized by a smooth
manifold~$M$, but instead for simplicity we work universally on~$\M X$.  As
in~\eqref{eq:41} fix a classifying map
  \begin{equation}\label{eq:61}
     f\:\M X\longrightarrow (MTSO_2,\Sigma \inv MTSO_1)\wedge  BG; 
  \end{equation}
then there is an induced classifying map 
  \begin{equation*}
     \fdot=r\circ f\:\M{\bX}\longrightarrow MTSO_1\wedge  BG 
  \end{equation*}
for $r$~the map~\eqref{eq:52}.  Set $\tau =\fdot^*(\lambda )$.  Let $t\:\M
X\to\M{\bX}$ be the restriction map on flat connections.  According to
Lemma~\ref{thm:6} the composition~$\tau \circ \sbuniv'\circ f$~is the reduced
twisting $\tau _t - (\dim\M X - t^*\dim\M{\bX})$.  The homotopy which
expresses commutativity of the square in~\eqref{eq:39} gives an isomorphism
  \begin{equation}\label{eq:43}
     \tau _t - (\dim\M X - t^*\dim\M{\bX})\;\overset{\cong }{\longrightarrow
     }\; -t^*\tau .  
  \end{equation}
In principle, $\dim\M X$ and $\dim\M{\bX}$~are locally constant functions
which vary over the moduli space.  However, the Euler characteristic of the
deformation complex~\eqref{eq:24} is independent of the connection, and we
easily deduce 
  \begin{equation*}
     \dim\M X - t^*\dim\M{\bX} \equiv (\dim G)\,b_0(\bX)\pmod2, 
  \end{equation*}
where $b_0$~is the number of connected components.  (We only track degrees in
$K$-theory modulo two; see~\eqref{eq:18}.)  According to the discussion
in~\S\ref{sec:2} (especially~\eqref{eq:68}), there then is an induced
pushforward
  \begin{equation}\label{eq:44}
     t_*\: K^{0}(\M X)\longrightarrow K^{\tau \,+\, (\dim G)\,b_0(\bX)}(\M{\bX}).  
  \end{equation}
This is the pushforward~\eqref{eq:8} associated to the bordism $X\:Y_0\to
Y_1$ with~$Y_0=\emptyset $ and~$Y_1=\bX$.  The invariant~\eqref{eq:8} is
then~$t_*(1)$. 
 
Note the special case~$\bX=\cir$.  Then $\M \cir=G\gpd G$ is the global
quotient stack of the action of~$G$ on~$G$ by conjugation.  The codomain
of~\eqref{eq:44} is thus
  \begin{equation*}
     K^{\tau +\dim G}(\M\cir)\cong K^{\tau +\dim G}(G\gpd G) = K^{\tau +\dim
     G}_G(G). 
  \end{equation*}
This is the basic space of the two-dimensional TQFT we construct;
see~\S\ref{sec:1}. 

Observe that the universal orientation, in the form described
around~\eqref{eq:39}, leads to the twist~$\tau $ in the codomain
of~\eqref{eq:44}.  This is the mechanism which was envisioned
in~\eqref{eq:68} when we discussed twistings in general: we have constructed
a pushforward from untwisted $K$-theory to twisted $K$-theory.  The universal
orientation neatly accounts for the construction of a Frobenius ring
structure on twisted $K$-theory.
 
To treat an arbitrary bordism $X\:Y_0\to Y_1$ we note that the deformation
complex at a flat connection~$a$ on a principal $G$-bundle $Q\to Y_0$ is 
  \begin{equation}\label{eq:54}
     \xymatrix{0 \ar[r] &\Omega ^0_{Y_0}(\mathfrak{g}\mstrut _{Q})
     \ar[r]^-{d_a} & \Omega ^1_{Y_0}(\mathfrak{g}\mstrut _{Q}) \ar[r] &0} 
  \end{equation}
The operator~$d_a$ is skew-adjoint.  Therefore, there is a canonical
trivialization of the $K$-theory class of~\eqref{eq:54}; e.g., a canonical
isomorphism $\ker d_a\cong \coker d_a$.  Suppose a classifying map~$f$ is
given as in~\eqref{eq:61} and let $\fdot_0,\fdot_1$ denote its restriction to
the boundary connections on~$Y_0,Y_1$.  Set $\tau _i=\fdot_i^*(\lambda )$.
Then \eqref{eq:43}~and the canonical trivialization of~\eqref{eq:54} on the
incoming boundary lead to the desired push-pull map 
  \begin{equation}\label{eq:71}
     \xymatrix{K^{\tau_0\,+\, (\dim G)\,b_0(Y_0)}(\M{Y_0}) \ar[r]^-{s^*} &
     K^{s^*\tau_0\,+\, (\dim G)\,b_0(Y_0)}(\M X) \ar[r]^-{t_*} &
     K^{\tau_1\,+\, (\dim G)\,b_0(Y_1)}(\M{Y_1})} 
  \end{equation}
This is the map~\eqref{eq:8} with the twistings induced from the universal
orientation.  
 
A universal orientation induces consistent orientations on the outgoing
restriction maps of bordisms.  In other words, if $X\:Y_0\to Y_1$ and
$X'\:Y_1\to Y_2$ are bordisms, then the push-pull maps derived from the
diagram  
  \begin{equation*}
     \xymatrix@!C{&&\M{X'\circ X}\ar[dl]_r \ar[dr]^{r'} \\
     &\M{X}\ar[dl]_s \ar[dr]^t&&\M{X'}\ar[dl]_{s'} \ar[dr]^{t'} \\
     \M{Y_0} && \M{Y_1} && \M{Y_2}} 
  \end{equation*}
satisfy 
  \begin{equation*}
     (t'r')_*\circ (sr)^* = [t'_*\circ {s'}^*]\circ [t_*\circ s^*]. 
  \end{equation*}
This follows from Lemma~\ref{thm:1} and the ``Fubini property'' 
  \begin{equation}\label{eq:59}
     (t'r')_* = t'_* r'_* 
  \end{equation}
of pushforward.  The orientation of~$t$ induces orientations of~$r'$
and~$t'r'$, since the diamond is a fiber product.  At stake in~\eqref{eq:59}
is the consistency of the orientations, which is ensured by the use of a
universal orientation.  The details of this argument\footnote{The main point
is that consistent orientations themselves form an \emph{invertible}
topological field theory, and these field theories factor through the group
completion of bordism, i.e., the Madsen-Tillmann space.} will be given on
another occasion.

One caveat: since $\M X,\M{\bX}$~are stacks we can only pushforward along
\emph{representable} maps, and this forces every component of~$X$ to have a
nonempty outgoing boundary.  As mentioned at the end of~\S\ref{sec:1}, the
partial topological quantum field theory obtained from the push-pull
construction extends to a full theory using the invertibility of the
(co)pairing attached to the cylinder~\cite[\S17]{FHT3}.

  \section{Families of Surfaces, Twistings, and Anomalies}\label{sec:4}

We begin with a general discussion of topological quantum field theories for
families.  Let $F$~be an $n$-dimensional TQFT in the most naive sense.  Thus
$F$~assigns a finite dimensional complex vector space~$F(Y)$ to a closed
oriented $(n-1)$-manifold~$Y$ and a linear map~$F(X)\:F(Y_0)\to F(Y_1)$ to a
bordism $X\:Y_0\to Y_1$.  In particular, $F(X)\in \CC$ if $X$~is closed.
Suppose that $\Y\to S$ is a fiber bundle with fiber a closed oriented
$(n-1)$-manifold.  Then the vector spaces assigned to each fiber fit together
into a complex vector bundle $F(\Y/S)\to S$.  If $\gamma \:[0,1]\to S$ is a
path, then $\gamma ^*\Y\to[0,1]$ is a bordism from the fiber~$\Y_{\gamma
(0)}$ to the fiber~$\Y_{\gamma (1)}$.  The topological invariance of~$F$
shows that $F(\gamma ^*\Y)\:F(\Y_{\gamma (0)}) \to F(\Y_{\gamma (1)})$ is
unchanged by a homotopy of~$\gamma $, and so $F(\Y/S)\to S$ carries a natural
flat connection.  Then a family of bordisms $\X\to S$ from $\Y_0\to S$ to
$\Y_1\to S$ produces a section~$F(\X/S)$ of $\Hom\bigl(F(\Y_0/S),F(\Y_1/S)
\bigr)\to S$; the topological invariance and gluing law of the TQFT imply
that this section is flat.  In other words, $F(\X/S)\in
H^0\bigl(S;\Hom\bigl(F(\Y_0/S),F(\Y_1/S) \bigr) \bigr)$.  It is natural,
then, to postulate that a TQFT in families gives more, namely classes of all
degrees: 
  \begin{equation}\label{eq:73}
     F(\X/S)\in H^{\bullet }\bigl(S;\Hom\bigl(F(\Y_0/S),F(\Y_1/S) \bigr)
     \bigr). 
  \end{equation}
These classes are required to satisfy gluing laws and topological invariance
as well as naturality under base change. 
 
The idea of a TQFT in families---at least in two dimensions---was introduced
in the mid~90s.  In two dimensions it is often formulated in a holomorphic
language (e.g.~\cite{KM}), and classes are required to extend to the
Deligne-Mumford compactification of the moduli space of Riemann surfaces.

Our push-pull construction works for families of surfaces---with a twist.
The purpose of this section is to alert the reader to the twist.\footnote{We
thank Veronique Godin for the perspicacious sign question which prompted this
exposition.}  
 
Suppose $\X\to S$ is a family of bordisms from $\Y_0\to S$ to $\Y_1\to S$,
where $\Y_i\to S$ are fiber bundles of oriented 1-manifolds.  Then the moduli
stacks of flat connections form a correspondence diagram over~$S$: 
  \begin{equation*}
     \xymatrix@!C{&\M{\X/S}\ar[dl]_s \ar[dr]^t \\ \M{\Y_0/S} \ar[dr]_{\pi _0}&&
     \M{\Y_1/S}\ar[dl]^{\pi _1} \\ &S} 
  \end{equation*}
The push-pull constructs a map from twisted~$K(\M{\Y_0/S})$ to
twisted~$K(\M{\Y_1/S})$.  We can also work locally over~$S$: the $K$-theory
of the fibers of~$\pi _i$ form bundles of spectra over~$S$ and the push-pull
gives a map of these spectra.  But for our purposes the global push-pull
suffices.  This construction is a variation of~\eqref{eq:73}: we use
$K$-theory rather than cohomology.
 
The discussion of~\S\ref{sec:3} goes through with one important change.  It
comes in the paragraph preceding~\eqref{eq:31}.  For simplicity suppose
$\Y_0=\emptyset $ so that the boundary~$\partial \X=\Y_1$ is entirely
outgoing.  Fix~$A\in \M{\X_s}$ a flat connection on a principal $G$-bundle
$P\to \X_s$.  Then the $KO$-theory class of the de Rham complex coupled to
the reduced adjoint bundle $\bg_P=\mathfrak{g}_P - \triv{\dim G}$ computes
the difference $T_A\M{\X_s} - (\dim G)\,H^{\bullet }(\X_s)$, where
$H^{\bullet }(\X_s)$~is the real cohomology of~$\X_s$ viewed as a class in
$KO$-theory.  In~\S\ref{sec:3} we treat~$H^{\bullet }(\X_s)$ as a trivial
vector space (there~$S=\pt$), but now $H^{\bullet }(\X_s)$~varies with~$s\in
S$ and so can give rise to a nontrivial twisting.  More precisely,
$H^{\bullet }(\X_s)$ is the fiber at~$s\in S$ of a flat vector bundle
$\mathscr{H}^{\bullet }(\X/S)\to S$.  Let $\tau _{\X/S}$~denote the twisting
of complex $K$-theory attached to this bundle.  This twisting replaces the
degree shift in~\eqref{eq:43} and the pushforward~\eqref{eq:44} is modified
to include that extra twist:
  \begin{equation}\label{eq:65}
     t_*\: K^{0}(\M {\X/S})\longrightarrow K^{\tau \,+\, (\dim G)\,\pi
     _1^*\tau_{\X/S}}(\M{\Y_1/S}).   
  \end{equation}
The degree shift is now incorporated into the twist~$\tau _{\X/S}$, and there
may be nontrivial contributions to the twist from~$w_1$ and~$W_3$ of
$\mathscr{H}^{\bullet }(\X/S)\to S$ as well.  (The degree shift and twistings
vanish canonically if $\dim G$~is even.)

        \begin{example}[]\label{thm:8}
 Consider the disjoint union~$X$ of two 2-disks.  The boundary circles are
outgoing, as above.  Suppose that $\dim G$~is odd.  For a single disk the
pushforward~$t_*(1)$ in~\eqref{eq:44} lands in $K^{\tau +1}_G(G)$ and is the
unit~$\unit$ in the Verlinde ring.  Thus for the disjoint union of two disks,
$t_*(1)$~is the image of~$\unit\otimes \unit$ under the external product
$K^{\tau +1}_G(G)\otimes K^{\tau +1}_G(G)\to K^{(\tau ,\tau )}_{G\times
G}(G\times G)$.  Now consider the family $\X\to S$ with fiber~$X$ and
base~$S=\cir$ in which the monodromy exchanges the two disks.  The flat
bundle $\mathscr{H}^{\bullet }(\X/S)\to S$ has rank~2 and nontrivial~$w_1$.
According to~\eqref{eq:65}, then, $t_*(1)$~for the family lives in the
twisted group $K^{(\tau ,\tau ) + \pi ^*\tau _{\X/S}}(\M{\partial \X/S})$.
On each fiber of $\pi \:\M{\partial \X/S}\to S$ we recover the
class~$\unit\otimes \unit$ above.  But upon circling the loop~$S=\cir$ this
class changes sign in the $\pi ^*\tau _{\X/S}$-twisted $K$-group.  Said
differently, the diffeomorphism which exchanges the disks acts by a sign
on~$\unit\otimes \unit$.  Of course, one might predict this from the sign
rule in graded algebra: the Verlinde ring~$K^{\tau +1}_G(G)$ is in odd
degree, so upon exchanging the factors of~$\unit\otimes \unit$ one picks up a
sign.  It shows up here as an extra twisting.
        \end{example}

This extra twisting is a topological analog of what is usually called an
\emph{anomaly} in quantum field theory.  In an anomalous theory in
$n$~dimensions the partition function on a closed $n$-manifold, rather than
being a complex-valued function on a space of fields, is a section of a
complex line bundle over that space of fields.  Furthermore, there is a gerbe
over the space of fields on a closed $(n-1)$-manifold, and for a bordism the
partition function is a map of the gerbes attached to the boundary.  In the
homological version described at the beginning of this section, the parameter
space~$S$ plays the role of the space of fields and for a family of closed
$n$-manifolds the partition function in a non-anomalous theory is an element
of~$H^{\bullet }(S;\RR)$.  An anomalous theory would assign a flat complex
line bundle $L\to S$ to the family, and the partition function would live in
the twisted cohomology $H^{\bullet }(S;L)= H^{L+\bullet }(S)$.  In the
2-dimensional TQFT we construct using push-pull on $K$-theory, the extra
$K$-theory twist~$\tau _{\X/S}$ is the anomaly; see~\eqref{eq:65}.  Notice
that there is no gerbe attached to a family of 1-manifolds (better: it is
canonically trivial).  We remark that the anomaly is itself a particular
example of an invertible topological quantum field theory.

\newpage 

 \bigskip\bigskip  

\providecommand{\bysame}{\leavevmode\hbox to3em{\hrulefill}\thinspace}
\providecommand{\MR}{\relax\ifhmode\unskip\space\fi MR }
\providecommand{\MRhref}[2]{%
  \href{http://www.ams.org/mathscinet-getitem?mr=#1}{#2}
}
\providecommand{\href}[2]{#2}

\end{document}